\documentclass[twocolumn]{autart}    % Enable this line and disable the
\usepackage{graphicx}
\graphicspath{{Figures/}}

\usepackage[utf8]{inputenc}
\usepackage[T1]{fontenc}
\usepackage[english]{babel}

\usepackage[retainorgcmds]{IEEEtrantools}
\usepackage{amssymb,amsmath}%,amsthm}
\usepackage{eucal}

\providecommand{\abs}[1]{\left\lvert#1\right\rvert}
\providecommand{\norm}[1]{\left\lVert#1\right\rVert}

\newcommand{\R}{\mathbb R}
\newcommand{\C}{\mathbb C}
\newcommand{\N}{\mathbb N}

\newcommand{\be}{\begin{equation}}
\newcommand{\ee}{\end{equation}}

\newcommand{\bes}{\begin{equation*}}
\newcommand{\ees}{\end{equation*}}

\newcommand{\ba}{\begin{eqnarray}}
\newcommand{\ea}{\end{eqnarray}}

\newcommand{\oi}{\overline{\imath}}
\newcommand{\oj}{\overline{\jmath}}

\newcommand{\on}{\overline{n}}
\begin{document}

\begin{frontmatter}
\runtitle{Null controllability of the heat equation}  % Running title for regular
                                              % papers but only if the title
                                              % is over 5 words. Running title
                                              % is not shown in output.

\title{Null controllability of the heat equation using flatness} % Title, preferably not more
                                                % than 10 words.

%\thanks[footnoteinfo]{This paper was not presented at any IFAC
%meeting. Corresponding author M.~T.~Cicero. Tel. +XXXIX-VI-mmmxxi.
%Fax +XXXIX-VI-mmmxxv.}

\author[First]{Philippe Martin}\ead{philippe.martin@mines-paristech.fr},  % Add the e-mail address (ead) as shown
\author[Second]{Lionel Rosier}\ead{Lionel.Rosier@univ-lorraine.fr},
\author[First]{Pierre Rouchon}\ead{pierre.rouchon@mines-paristech.fr}

\address[First]{Centre Automatique et Systèmes, MINES ParisTech, 60 boulevard Saint-Michel
75272 Paris Cedex 06, France}
\address[Second]{Institut Élie Cartan de Lorraine, UMR 7502 UdL/CNRS/INRIA, BP 70239, 54506 Vandœuvre-lès-Nancy, France}


\begin{abstract}                          % Abstract of not more than 200 words.
We derive in a direct and rather straightforward way the null controllability of the $N$-dimensional heat equation in a bounded cylinder with boundary control at one end of the cylinder. We use the so-called \emph{flatness approach}, which consists in parameterizing the solution and the control by the derivatives of a ``flat output''. This yields an explicit control law achieving the exact steering to zero. Replacing the involved series by partial sums we obtain a simple numerical scheme for which we give explicit error bounds. Numerical experiments  demonstrate the relevance of the approach.
\end{abstract}

\end{frontmatter}
%===============================================================================

\section{Introduction}
The controllability of the heat equation was first considered in the 1-D case~\cite{FattoR1971ARMA,LuxemK1971TAMS,Jones1977JMAA,Littm1978ASNSPCS} and very precise results were obtained by the classical moment approach.  Next using Carleman estimates and duality arguments the null controllability was proved in \cite{FursiI1996book,LebeaR1995CPDE} for any bounded domain in~$\R^N$, any control time~$T$, and any control region.
This Carleman approach proves very efficient also with semilinear parabolic equations~\cite{FursiI1996book}.
%This Carleman approach has also proved very efficient to deal with semilinear parabolic equations, e.g. the dissipative Burger equation, the Navier-Stokes equation~\cite{FursiI1996book}, the Ginzburg-Landau equation~\cite{RosieZ2009AIHPAANL}; see also \cite{AmmarBGT2011MCRF} for a survey of null controllability results for systems of parabolic equations.
%
%
By contrast the interest for the numerical investigation of the null controllability of the heat equation (or of parabolic equations) is fairly recent: apart from~\cite{CarthGL1994JOTA}, the first significant contributions are~\cite{Zuazu2006,Zheng2008AA,MunchZ2010IP,BoyerHL2011NM,MicuZ2011SCL,BelgaK2011IP,FernaM2011HAL}; see also~\cite{GarciOTToappearJoIaIP} for an application to some inverse problems.
%%
%A natural candidate for the control input is the control of minimal $L^2-$norm, which may be obtained as a trace of the solution of the (backward) adjoint problem whose terminal state is the minimizer of a suitable quadratic cost. Unfortunately its computation is a hard task~\cite{MicuZ2011SCL}; indeed the terminal state of the adjoint problem associated with some regular initial state of the control problem may be highly irregular, which leads to severe troubles in the numerical computation of the control function.
%
All the above results rely on some observability inequalities for the adjoint system. A direct approach which does not involve the adjoint problem was proposed in~\cite{Jones1977JMAA,Littm1978ASNSPCS,LinL1995AMO,LittmT2007}. In~\cite{Jones1977JMAA} a fundamental solution for the heat equation with compact support in time was introduced and used to prove null controllability.
%The same approach was used in~\cite{Rosie2002CAM} to derive a fundamental solution with compact support in time for some dispersive equations (including the Schrödinger equation and the Korteweg-de Vries equations), and used to prove some controllability results in unbounded domains.
The results in \cite{Jones1977JMAA,Rosie2002CAM} can be used to derive control results on a bounded interval with two or one boundary control in some Gevrey class, or on a bounded domain of~$\R^N$ with a control supported on the whole boundary (see also~\cite{LittmT2007}). An extension of those results to the semilinear heat equation in 1-D was obtained in~\cite{LinL1995AMO} in a more explicit way through the resolution of an ill-posed problem with data of Gevrey order~2 in~$t$.

In this paper, which builds on the preliminary versions~\cite{MartiRR2013CPDE,MartiRR2013CDC}, we derive in a straightforward way the null controllability of the heat equation in a bounded cylinder $\Omega=\omega\times(0,1)\subset\R^N$ with Neumann boundary control on~$\omega\times\{1\}$. More precisely given any final time~$T>0$ and initial state $\theta_0\in L^2(\Omega)$ we provide an explicit and very regular control such that the state reached at time~$T$ is exactly zero.
We use the so-called \emph{flatness} approach~\cite{FliesLMR1995IJoC}, which consists in parameterizing the solution~$\theta$ and the control~$u$ by the derivatives of a ``flat output''~$y$; this notion was initially introduced for finite-dimensional (nonlinear) systems, and later extended to (in particular) parabolic 1-dimensional PDEs~\cite{LarocMR2000IJRNC,LynchR2002IJC,MeureZ2008MCMDS,Meure2011A}. Choosing a suitable trajectory for this flat output~$y$ then yields an explicit series for a control achieving the exact steering to zero.
Note this paper is probably the first example of using flatness for the motion planning of a ``truly'' $N$-dimensional PDE.
Comparing our results to~\cite{LinL1995AMO,LittmT2007} we note that: (i)~our control is not supported on the whole boundary even in dimension~$N>1$; (ii)~the control and the solution are Gevrey of order $s\in (1,2)$ in time; (iii)~the control and the solution are developed in series whose easily computed partial sums yield accurate numerical approximations of both the control and the solution.

The paper runs as follows. In section~\ref{sec:dim1} we consider the control problem in dimension~$N=1$. In Proposition \ref{prop1} we investigate an ill-posed problem with Cauchy data in a Gevrey class and prove its (global) well-posedness. Theorem~\ref{thm1} then establishes the null controllability in small time for any initial data in~$L^2$.
In section~\ref{sec:dimN} we extend these results to the cylinder $\Omega=\omega\times(0,1)\subset\R^N$.
Section~\ref{sec:estimates} provides accurate error estimates when the various series involved are replaced by their partial sums. Finally in section~\ref{sec:numerics} some numerical experiments demonstrate the interest of the approach.

\section{Preliminaries (Gevrey functions)}\label{eq:preliminaries}
In the sequel we consider series with infinitely many derivatives of functions. The notion of Gevrey order is a way of estimating the growth of these derivatives: we say that a function $y\in C^\infty([0,T])$ is {\em Gevrey of order $s\geq0$ on~$[0,T]$} if there exist positive constants~$M,R$ such that
\bes
%\label{G1}
\abs{y^{(p)}(t)} \leq M\frac{p!^s}{R^p} \qquad \forall t\in [0,T],\ \forall p\ge 0.
\ees
More generally if $K\subset\R^N$ is a compact set and $y$ is a function of class~$C^\infty$ on~$K$ (i.e. $y$ is the restriction to $K$ of a function of class~$C^\infty$ on some open neighbourhood $\Omega$ of~$K$), we say $y$ is {\em Gevrey of order $s_1$ in $x_1$, $s_2$ in $x_2$,\ldots,$s_N$ in $x_N$ on~$K$} if there exist positive constants $M,R_1,...,R_N$ such that $\forall x\in K,\ \forall p\in\N^N$
\bes
%\label{G3}
\abs{\partial_{x_1}^{p_1}\partial_{x_2}^{p_2}\cdots\partial_{x_N}^{p_N}y(x)} \le
M\frac{\prod_{i=1}^N (p_i !)^{s_i}}{\prod_{i=1}^N R_i^{p_i} }.%, \quad \forall x\in K,\ \forall p\in\N^N.
\ees
By definition, a Gevrey function of order $s$ is also of order $r$ for~$r\geq s$. Gevrey functions of order~1 are analytic (entire if $s<1$).  Gevrey functions of order~$s>1$ have a divergent Taylor expansion; the larger~$s$, the ``more divergent'' the Taylor expansion. Important properties of analytic functions generalize to Gevrey functions of order $s>1$: the scaling, addition, multiplication and derivation of Gevrey functions of order $s>1$ is of order~$s$, see~\cite{Ramis1978,Rudin1987book,Yaman1989AGAG}. But contrary to analytic functions, Gevrey functions of order $s>1$ may be constant on an open set without being constant everywhere.  For example the ``step function''
\bes
\phi_s(t):=\begin{cases}
1 & \text{if $t\leq0$}\\
0 & \text{if $t\geq1$}\\
\dfrac{ e^{-(1-t)^{-k}} }{ e^{-(1-t)^{-k}} + e^{-t^{-k}} }
&\text{if $t\in(0,1)$},
\end{cases}
\ees
where $k = (s-1)^{-1}$ is Gevrey of order~$s$ on $[0,1]$ (and in fact on~$\R$); notice $\phi_s(0)=1$, $\phi_s(1)=0$ and $\phi_s^{(i)}(0)=\phi_s^{(i)}(1)=0$ for all~$i\geq1$.

In conjunction with growth estimates we will repeatedly use Stirling's formula $n!\sim(n/e)^n\sqrt{2\pi n}$.

\section{The one-dimensional heat equation}\label{sec:dim1}
For simplicity we first study the 1-D heat equation
\begin{IEEEeqnarray}{rCl'l}
    \partial_t\theta(t,x)-\partial^2_{x}\theta(t,x) &=& 0, &(t,x)\in(0,T)\times(0,1)\quad \label{B1}\\
    \partial_x\theta(t,0) &=& 0, &t\in(0,T) \label{B2}\\
    \partial_x\theta(t,1) &=& u(t), &t\in(0,T) \label{B3}
\end{IEEEeqnarray}
with initial condition in~$L^2(0,1)$
\begin{IEEEeqnarray*}{rCl'l}
    \theta(0,x) &=&\theta_0(x), &x\in (0,1). %\label{B4}
\end{IEEEeqnarray*}

We claim the system~\eqref{B1}--\eqref{B3} is ``flat'' with $y(t):=\theta(t,0)$ as a flat output, meaning there is (in appropriate spaces of smooth functions) a $1-1$ correspondence between arbitrary functions $t\mapsto y(t)$ and solutions of~\eqref{B1}--\eqref{B3}.

We first seek a formal solution in the form
\bes
\theta(t,x):=\sum_{i\ge0}\frac{x^i}{i!}a_i(t)
\ees
where the $a_i$'s are functions yet to define. Plugging this expression into~\eqref{B1} yields
\bes
\sum_{i\ge 0}\frac{x^i}{i!}[a_{i+2}-a_i']=0,
\ees
hence $a_{i+2}=a_i'$ for all~$i\ge0$. On the other hand $y(t)=\theta(t,0)=a_0(t)$, and \eqref{B2} implies~$a_1(t)=0$. As a consequence $a_{2i}=y^{(i)}$ and $a_{2i+1}=0$ for all~$i\ge0$. The formal solution thus reads
\be
\label{AA10}
\theta(t,x)=\sum_{i\ge 0}\frac{x^{2i}}{(2i)!}y^{(i)}(t)
\ee
while the formal control is given by
\be
\label{AA10bis}
u(t)=\theta_x(t,1)=\sum_{i\ge 1}\frac{y^{(i)}(t)}{(2i-1)!}.
\ee

We now  give a meaning to this formal solution by restricting $t\mapsto~y(t)$
to be Gevrey of order~$s\in [0,2)$.
\begin{prop}
\label{prop1}
Let $s\in [0,2)$, $-\infty <t_1<t_2<\infty$, and $y\in C^\infty ([t_1,t_2] )$ satisfying for some constants $M,R>0$
\be
\label{AA11}
\abs{y^{(i)}(t)} \le M \frac{i!^s}{R^i}, \qquad \forall i\ge 0,\ \forall t\in [t_1,t_2].
\ee
Then the function $\theta$ defined by \eqref{AA10} is Gevrey of order $s$ in $t$ and $s/2$ in $x$ on
$[t_1,t_2]\times [0,1]$; hence the control~$u$ defined by \eqref{AA10bis} is also Gevrey of order $s$ on $[t_1,t_2]$.
\end{prop}
\begin{pf}
We must prove the formal series
\be
\label{AA12}
\partial _t ^m\partial _x ^n \theta (t,x) = \sum_{2i\ge n}  \frac{ x^{2i-n } }{ (2i-n)! } y^{(i+m)}(t)
\ee
is uniformly convergent on $[t_1,t_2]\times [0,1]$ with growth estimates of the form
\be
\label{AA12bis}
\abs{\partial _t ^m\partial _x ^n \theta (t,x)} \le C
\frac{m! ^s}{R_1^m}\, \frac{n! ^\frac{s}{2}}{R_2^n} \cdot
\ee
By \eqref{AA11}, we have for all $(t,x)\in [t_1,t_2]\times [0,1]$
\begin{eqnarray*}
\left\vert  \frac{x^{2i-n}}{(2i-n)!} y^{(i+m)} (t) \right\vert &\le& \frac{M}{R^{i+m}} \, \frac{(i+m)! ^s}{(2i-n)!} \\
&\le& \frac{M}{R^{i+m} }\,  \frac{(2^{i+m} i! \, m! )^s}{(2i-n)!} \\
&\le& \frac{M'}{R^{i+m}}\, \frac{2^{si} \bigl(2^{-2i} \sqrt{\pi i} \, (2i)! \bigr)^\frac{s}{2} }{(2i-n)!}\, \frac{m!^s}{2^{-sm}} \\
&\le&  M' \frac{(\pi i) ^{\frac{s}{4}} }{R_1 ^i (2i-n)! ^{1-\frac{s}{2}} } n! ^{\frac{s}{2}} \frac{m!^s}{R_1 ^m},
\end{eqnarray*}
where we have set $R_1 = 2^{-s}R$ and $M'\geq M$ is a suitable constant; we have used Stirling's formula for $(2i)!$ and twice $(i+j)! \le 2^{i+j} i!j!$. Since $\sum_{2i\ge n} \frac{(\pi  i)^\frac{s}{4} }{R_1 ^i (2i-n)! ^{1-\frac{s}{2} }  } <\infty$ the series in~\eqref{AA12} are uniformly convergent for all $m,n\ge 0$, hence $\theta \in C^\infty ([t_1,t_2]\times [0,1])$.
Finally, since
\bes
\sum_{2i\ge n} \frac{ M'(\pi  i)^\frac{s}{4} }{R_1 ^i (2i-n)! ^{1-\frac{s}{2} }  }
\le  \frac{M'}{R_1^{\frac{n}{2}}}\Bigl(\frac{\pi}{2}\Bigr)^\frac{s}{4}  \sum_{j\ge 0} \frac{ j^\frac{s}{4} + n^\frac{s}{4} }{R_1 ^\frac{j}{2}  j! ^{1-\frac{s}{2} }  }
\le C R_2^{-n}
\ees
where $R_2\in (0,\sqrt{R_1})$ and $C>0$ is some constant independent of~$n$, we have the desired estimates~\eqref{AA12bis}.\qed
\end{pf}

\begin{rem}
Note the Cauchy-Kovalevsky theorem, see e.g.~\cite{Horma1983book}, only ensures for $y\in C^\omega(0,T)$ the existence of solutions of the Cauchy problem composed of~\eqref{B1}-\eqref{B2} and $\theta(t,0)=y(t)$ on a small neighborhood of~$\{x=0\}$; whereas we are interested in non-analytic solutions~$\theta$ defined for all $(t,x)\in[0,T]\times [0,1]$.
\end{rem}

%==============================================================================================
We now derive an explicit control steering the system from any initial state~$\theta_0\in L^2(0,1)$ at time~$0$ to the final state~$0$ at time~$T>0$. Two ideas are involved: on the one hand thanks to the flatness property it is easy to find a control achieving the steering to zero starting from a certain set of initial conditions (step~2 in the proof of Theorem~\ref{thm1}); on the other hand thanks to the regularizing property of the heat equation this set is reached from any~$\theta_0\in L^2(0,1)$ when applying first a zero control for some time (step~1 in the proof of Theorem~\ref{thm1}).

\begin{thm}\label{thm1}
Let $\theta_0\in L^2(0,1)$, $T>0$, $\tau\in(0,T)$ and $s\in(1,2)$. Then there exists a function $y:[\tau,T]\to\R$ Gevrey of order~$s$ on~$[\tau,T]$ such that the control
\bes
u(t):=\begin{cases}
0 & \text{if $0\le t\le\tau$}\\
\sum_{i\ge1}\frac{y^{(i)}(t)}{(2i-1)!} &\text{if $\tau <t\le T$}
\end{cases}
\ees
steers the system~\eqref{B1}--\eqref{B3} from the initial state $\theta(0,.)=\theta_0$ at time~$0$ to the final state~$\theta(T,.)=0$ at time~$T$.

Moreover $u$ is Gevrey of order~$s$ on~$[0,T]$;
%$t\mapsto\theta(t,\cdot)$ is in $C\bigl([0,T],L^2(0,1)\bigr)$;
the solution $\theta$ of~\eqref{B1}--\eqref{B3} is Gevrey of order~$s$ in~$t$ and~$s/2$ in~$x$ on~$[\varepsilon,T]\times[0,1]$ for all $\varepsilon\in(0,T)$.
\end{thm}
\begin{pf} We first apply a null control on~$[0,\tau ]$ to reach a regular intermediate state (step~1), and then use the flatness property on~$[\tau , T]$ to steer this intermediate state to~0 (step~2).
\paragraph*{\sc Step~1: free evolution}
Let $u(t)=0$ for $t\in [0,\tau ]$. Decompose $\theta _0$ as the Fourier series of cosines
\bes
\theta_0(x)=\sum_{n\ge 0}c_n\sqrt{2}\cos(n\pi x)
\ees
where the convergence holds in $L^2(0,1)$ and
\bes
2|c_0|^2+\sum_{n\ge 1} |c_n|^2 =\int_0^1 |\theta _0(x)|^2 dx <\infty.
\ees
The solution starting from~$\theta_0$ then reads
\be\label{C100}
\theta (t,x)=\sum_{n\ge 0}c_ne^{-n^2\pi^2t}\sqrt{2}\cos(n\pi x)
\ee
and in particular
\be\label{B10}
\theta_\tau(x):=\theta(\tau,x)=\sum_{n\ge 0}c_n e^{-n^2\pi^2\tau}\sqrt{2}\cos(n\pi x).
\ee

\begin{lem}\label{lem1}
$\theta_\tau$ is analytic in $\C$ and can be expanded as
\bes%\label{B11}
\theta_\tau(x) = \sum_{i\ge0}y_i\frac{x^{2i}}{(2i)!}, %\quad x\in\C,
\ees
with
\bes%\label{B12}
%\abs{y_i} \le C\frac{i!}{\tau^{i+\frac{1}{2}}}
\abs{y_i} \le C\Bigl(1+\frac{1}{\sqrt\tau}\Bigr)\frac{i!}{\tau^i}
\ees
and $0<C\le K\sup_{n\ge 0} |c_n|$, where $K$ is a universal positive constant.
\end{lem}
\noindent\textbf{Proof of Lemma~\ref{lem1}:~}
$\theta_\tau$ is analytic in~$\C$ since for~$\abs{x}\leq r$
\bes
\abs{c_ne^{-n^2\pi^2\tau}\sqrt{2}\cos(n\pi x)} \le Ce^{-n^2\pi^2\tau+n\pi r},
\ees
which ensures the uniform convergence of the series~\eqref{B10} in every open disk of radius $r>0$. On the other hand
\begin{IEEEeqnarray*}{rCl}
\theta_\tau(x) &=& \sqrt{2}\sum_{n\ge 0}c_ne^{-n^2\pi^2\tau}\sum_{i\ge0}(-1)^i\frac{(n\pi x)^{2i}}{(2i)!}\\
&=& \sum_{i\ge0} \frac{x^{2i}}{(2i)!} \underbrace{\left(\sqrt{2}\sum_{n\ge0}c_n e^{-n^2\pi^2\tau}(-n^2\pi^2)^i\right)}_{=:y_i}.
\end{IEEEeqnarray*}
The change in the order of summation will be justified once proved that
\bes
\sum_{i,n\ge0} \frac{(n\pi x)^{2i}}{(2i)!} e^{-n^2\pi^2\tau} <\infty,  \quad\forall x\geq0.
\ees
For $i\ge 0$ let $h_i(x):=e^{-\tau\pi^2x^2}(\pi x)^{2i}$ and $N_i:=\left[\bigl(\frac{i}{\pi ^2\tau}\bigr)^{\frac{1}{2}}\right]$. The map $h_i$ is increasing on $\bigl[0,\left(\frac{i}{\pi ^2\tau}\right)^{\frac{1}{2}}\bigr]$ and decreasing on $\bigl[\left(\frac{i}{\pi ^2\tau}\right)^{\frac{1}{2}},+\infty\bigr)$ hence
\begin{IEEEeqnarray*}{rCl}
\sum_{n\ge0}h_i(n) &\le& \int_0^{N_i}h_i(x)dx + h_i(N_i)\\
&& +\: h_i (N_i + 1) + \int_{N_{i+1} }^\infty h_i(x)dx\\
&\le& 2h_i\left(\Bigl(\frac{i}{\pi^2\tau}\Bigr)^{\frac{1}{2}}\right) + \int_0^\infty h_i(x)dx\\
&\le& C\frac{i!}{\tau^i} + \int_0^\infty h_i(x)dx;
\end{IEEEeqnarray*}
we have used Stirling's formula for the last line. On the other hand integrating by parts yields
\begin{IEEEeqnarray*}{rCl}
\int_0^\infty h_i(x)dx &=& \frac{2i-1}{2\tau}\int_0^\infty h_{i-1}(x)dx\\
&=& \frac{(2i-1)\cdots 3\cdot1}{(2\tau)^i}\int_0^\infty e^{-\tau\pi^2x^2}dx\\
&=& \frac{(2i)!}{2^ii!(2\tau)^i}\cdot\frac{1}{\pi\sqrt\tau} \int_0^\infty e^{-x^2}dx\\
&\leq& C\frac{i!}{\tau^i\sqrt{\tau}},
\end{IEEEeqnarray*}
where we have again used Stirling's formula. Therefore
\bes
\abs{y_i}\le\sqrt{2}\sup_{n\ge 0}\abs{c_n}\sum_{n\ge0}h_i(n) \le C\Bigl(1+\frac{1}{\sqrt\tau}\Bigr)\frac{i!}{\tau^i}
\ees
Moreover
\bes
\sum_{i,n\ge0} \frac{(n\pi x)^{2i}}{(2i)!} e^{-n^2\pi^2\tau}
\leq C\Bigl(1+\frac{1}{\sqrt\tau}\Bigr)\sum_{i\ge0}\underbrace{\frac{i!}{(2i)!}\Bigl(\frac{x^2}{\tau}\Bigr)^i}_{=:v_i}<\infty
\ees
since $\frac{v_{i+1}}{v_i}\sim\frac{1}{4i}\frac{x^2}{\tau}$. Lemma~\ref{lem1} is proved.\qed

The following lemma establishes another aspect of the regularizing properties of the heat equation.
\begin{lem}\label{lem:regstep1}
The solution $\theta$ defined by~\eqref{C100} is Gevrey of order~$1$ in~$t$ and $1/2$ in~$x$
on $[\varepsilon,\tau]\times [0,1]$ for all $\varepsilon\in(0,\tau)$.
\end{lem}
\noindent\textbf{Proof of Lemma~\ref{lem:regstep1}:~}
first notice that for all $\delta>0$
\begin{IEEEeqnarray*}{rCl}
\IEEEeqnarraymulticol{3}{l}{
(\delta k^2\pi^2)^{m+\frac{n}{2}}e^{-\frac{3\delta}{2}k^2\pi^2}
}\\ \quad
&=& \Bigl((\delta k^2\pi^2)^me^{-\delta k^2\pi^2}\Bigr)
\Bigl((\delta k^2\pi^2)^ne^{-\delta k^2\pi^2}\Bigr)^\frac{1}{2}
\\
&\leq& m!n!^{\frac{1}{2}}
\end{IEEEeqnarray*}
by using $x^p/p!\leq e^x$ twice. Hence for all $m,n\ge0$
%and $(t,x)\in[\varepsilon,\tau]\times[0,1]$
\begin{IEEEeqnarray*}{rCl}
\IEEEeqnarraymulticol{3}{l}{
\abs{\partial_t^m\partial_x^n\theta(t,x)}
}\\ \quad
&\leq&
\sqrt{2}\sup_{k\ge0}\abs{c_k}\sum_{k\ge 0}(k\pi)^{2m+n}e^{-\varepsilon k^2\pi^2}
\\
&\leq& \sqrt{2}\sup_{k\ge0}\abs{c_k}\sum_{k\ge 0}
\frac{(\delta k^2\pi^2)^{m+\frac{n}{2}}e^{-\frac{3\delta}{2}k^2\pi^2}}{\delta^{m+\frac{n}{2}}}
e^{-(\varepsilon-\frac{3\delta}{2})k^2\pi^2}
\\
&\leq& \biggl(\sqrt{2}\sup_{k\ge0}\abs{c_k}\sum_{k\ge 0}e^{-(\varepsilon-\frac{3\delta}{2})k^2\pi^2}\biggr)
\frac{m!}{\delta^m}\frac{n!^\frac{1}{2}}{\delta^\frac{n}{2}}.
\end{IEEEeqnarray*}
Picking any $\delta<\frac{2\varepsilon}{3}$ the last series is obviously convergent, which proves $\theta$ is Gevrey of order~$1$ in~$t$ and $1/2$ in~$x$ on~$[\varepsilon,\tau]\times [0,1]$ for all $\varepsilon\in(0,\tau)$.\qed

\paragraph*{\sc Step~2: construction of the control on~$[\tau,T]$}
\begin{lem}\label{lem2}
Let $0<\tau<T$ and $1< s< 2$. The function $y:[\tau,T]\to\R$
\bes
y(t):=\phi_s\Bigl(\frac{t-\tau}{T-\tau}\Bigr)\sqrt{2}\sum_{n\ge0}c_n e^{-n^2\pi^2t}.
\ees
is Gevrey of order $s$ on $[\tau,T]$ and satisfies for all $i\ge 0$
\begin{IEEEeqnarray}{rCl}
y^{(i)}(\tau) &=& y_i\label{Q1}\\
y^{(i)}(T) &=& 0\label{Q2}\\
\abs{y^{(i)}(t)} &\le&  C\sup_{n\ge0}\abs{c_n}\frac{i!^s}{R^i}, \quad\forall t\in[\tau,T] \label{Q3}
\end{IEEEeqnarray}
for some constants $C,R>0$ depending only on~$\tau,T,s$; $\phi_s$ is a Gevrey ``step function'' (see introduction) and the $y_i$'s are as in Lemma~\ref{lem1}.
\end{lem}
\noindent\textbf{Proof of Lemma~\ref{lem2}:~}
$\tilde y(t):=\sqrt{2}\sum_{n\ge0}c_n e^{-n^2\pi^2t}$ is clearly analytic on $(0,+\infty )$ hence Gevrey of order~$s$ on~$[\tau,T]$, and satisfies for all~$i\ge0$
\bes
\tilde y^{(i)}(\tau) =  \sqrt{2}\sum_{n\ge0}c_ne^{-n^2\pi^2\tau}(-n^2\pi^2)^i = y_i.
\ees
As a consequence $y$ is Gevrey of order $s$ on~$[\tau,T]$, and satisfies $y^{(i)}(\tau)=y_i$ and $y^{(i)}(T)=0$.

Finally for all $t\in [\tau,T]$
\begin{IEEEeqnarray*}{rCl}
\abs{{\bar y}^{(i)}(t)}
&\le& \sqrt{2}\sup_{n\ge 0}\abs{c_n}\sum_{n\ge 0} e^{-n^2\pi ^2 \tau} (n^2\pi ^2)^i\\
&\le& \sqrt{2} \sup_{n\ge 0}\abs{c_n}\frac{i!}{ (\frac{\tau}{2} )^i}
\sum_{n\ge0}e^{-n^2\pi^2\frac{\tau}{2}}, %\label{QQ1}
\end{IEEEeqnarray*}
where we have used $x^i/i!\le e^x$ for $x=n^2\pi^2\tau/2$; \eqref{Q3} then follows from~\cite[Theorem~$19.7$]{Rudin1987book}. The proof of Lemma \ref{lem2} is complete. \qed

The control defined on~$[\tau,T]$ by
\bes
u(t):=\sum_{i\ge1}\frac{y^{(i)}(t)}{(2i-1)!}
\ees
steers the system~\eqref{B1}--\eqref{B3} from $\theta_\tau$ at time~$\tau$ to 0 at time~$T$. Indeed the solution~$\theta$ on~$[\tau,T]$ is given by~\eqref{AA10}, and by~\eqref{Q1}-\eqref{Q2} satisfies $\theta(\tau,\cdot)=\theta_\tau$ and $\theta(T,\cdot)=0$. By Proposition~\ref{prop1} this solution $\theta$ is Gevrey of order $s$ in~$t$ and $s/2$ in~$x$ on~$[\tau,T]\times[0,1]$; on the other hand by Lemma~\ref{lem:regstep1} the solution $\theta$ defined by~\eqref{C100} is Gevrey of order~$1$ in~$t$ and $1/2$ in~$x$ on $[\varepsilon,\tau]\times [0,1]$ for all $\varepsilon\in(0,\tau)$.
Since $s>1$, to prove $\theta$ is Gevrey of order $s$ in~$t$ and $s/2$ in $x$ on $[\varepsilon,T]\times[0,1]$ for all $\varepsilon\in(0,T)$ it is then sufficient to check
$\partial_t^k\theta(\tau^+,x)=\partial_t^k\theta(\tau^-,x)$ for $k\geq0$ and~$x\in[0,1]$.
But by~\eqref{Q1} the series \eqref{AA10} and~\eqref{C100} coincide at $t=\tau$, hence so do their (space) derivatives. Therefore
\bes
\partial_t^k\theta(\tau^+,x)
=\partial^{2k}_x\theta(\tau^+,x)
=\partial^{2k}_x\theta(\tau^-,x)
=\partial_t^k\theta(\tau^-,x).
\ees
%%%%%%%%%%%%%%% OLD EXPLICIT COMPUTATION
%\begin{IEEEeqnarray*}{rCl}
%\IEEEeqnarraymulticol{3}{l}{
%\partial_t^k\theta(\tau^+,x)
%}\\ \quad
%& = & \sum_{i\ge 0}\frac{x^{2i}}{(2i)!}y^{(i+k)}(\tau)\\
%%
%&=& \sum_{i\ge 0} \frac{x^{2i}}{(2i)!}y_{i+k}\\
%%
%&=& \sqrt{2} \sum_{i\ge 0} \frac{x^{2i}}{(2i)!}
%\left( \sum_{n\ge 0} c_n e^{-n^2\pi ^2 \tau } n^{2(i+k)}  \right) (-\pi^2)^{i+k} \\
%%
%&=& \sum_{n\ge 0}  c_n (-n^2\pi ^2) ^k e^{-n^2\pi ^2 \tau }  \sqrt{2} \cos (n\pi x)\\
%%
%&=& \partial_t^k\theta(\tau^-,x).
%\end{IEEEeqnarray*}
As a simple consequence $u=\theta_x(t,1)$ is Gevrey of order~$s$ on~$[0,T]$, which concludes the proof of Theorem~\ref{thm1}. \qed
\end{pf}

\section{The N-dimensional heat equation}\label{sec:dimN}
Let $\omega\subset\R^{N-1}$, $N\ge2$, be a smooth\footnote{Most results in this section are valid also for $\omega =(0,1)^{N-1}$} bounded open set, and $\Omega:=\omega\times (0,1)\subset \R ^N$;
$x\in\Omega$ is written $x=(x',x_N)$ with $x'=(x_1,\dots,x_{N-1})\in \omega$ and $x_N\in (0,1)$.
We are concerned with the null controllability of the system
\begin{IEEEeqnarray}{rCl'l}
    \partial_t\theta(t,x) &=& \Delta\theta(t,x), &t\in(0,T),x\in\Omega \label{D1}\\
    \partial_\nu\theta(t,x',1) &=& u(t,x'), &t\in(0,T),x'\in\omega \label{D2}\\
    \partial_\nu\theta(t,x) &=& 0, &t\in(0,T),x\in\partial\Omega\setminus\omega\times\{1\}
    \IEEEeqnarraynumspace\label{D3}
\end{IEEEeqnarray}
($\Delta:=\partial_{x_1}^2+\cdots+\partial_{x_N}^2$, $\nu$ is the outward unit normal vector to~$\partial\Omega$), with initial condition in~$L^2(\Omega)$
\begin{IEEEeqnarray*}{rCl'l}
    \theta(0,x) &=&\theta_0(x), &x\in\Omega. %\label{D4}
\end{IEEEeqnarray*}

Let $(e_j)_{j\ge0}$ be an orthonormal basis of $L^2(\omega)$ such that each $e_j$ is an eigenfunction for the Neumann Laplacian on $\omega$. In other words $e_j\in H^2(\omega)$ and
\begin{IEEEeqnarray*}{rCl's}
-\Delta'e_j(x') &=& \lambda_je_j(x'), & $x'\in\omega$\\
\frac{\partial e_j}{\partial\nu'}(x') &=&0, & $x'\in\partial\omega$,
\end{IEEEeqnarray*}
where $\Delta':=\partial_{x_1}^2+\cdots+\partial_{x_{N-1}}^2$,
$\nu  '$ is the outward unit normal vector to $\partial \omega$,  $0=\lambda _0<\lambda _1\le \lambda _j\le \lambda _{j+1}$, and
$e_0(x)=|\omega |^{-1/2}$.
Note that by Weyl's formula, see e.g.~\cite[Chapter~8]{Roe1998book}\cite[Section~15]{Shubi2001book},
there are constants $A_1,A_2>0$ such that
\be\label{weyl}
A_1j^\frac{2}{N-1}\leq\lambda_j\leq A_2j^\frac{2}{N-1},\qquad j\ge0.
\ee

Now decompose $\theta(t,x',0)$ as
\be\label{D100}
    \theta(t,x',0) =\sum_{j\geq0} z_j(t)e_j(x').
\ee
We claim \eqref{D1}-\eqref{D3} is ``flat'' with $z(t):=\bigl (z_j(t)\bigr) _{j\geq0}$ as a flat output, i.e.
there is a 1-1 correspondence between arbitrary functions $t\mapsto z(t)$ in a certain space and
smooth solutions $\theta$ of~\eqref{D1}--\eqref{D3}. Indeed the $z_k$'s hence $z$ are obviously uniquely defined by a solution $\theta$ of~\eqref{D1} since
\begin{IEEEeqnarray*}{rCl}
    \int_{\omega}\theta(t,x',0)e_k(x')dx' &=&\sum_{j\geq0}z_j(t)\int_{\omega}e_j(x')e_k(x')dx'\\
    &=& z_k(t).
\end{IEEEeqnarray*}
Conversely, given a sequence $z(t)=(z_j(t))_{j\geq0}$ of functions in $C^\infty ([0,T])$ we
seek a formal solution in the form
\begin{IEEEeqnarray*}{rCl}
    \theta(t,x',x_N) &=&\sum_{i\geq0}\frac{x_N^i}{i!}a_i(t,x'),
\end{IEEEeqnarray*}
where the $a_i$'s are functions yet to be defined. Plugging this formal solution into~\eqref{D1} we find
\begin{IEEEeqnarray*}{rCl}
    \sum_{i\geq0}\frac{x_N^i}{i!}\bigl[a_{i+2}(t,x')-(\partial_t-\Delta ')a_i(t,x')\bigr] &=&0,
\end{IEEEeqnarray*}
hence $a_{i+2}=(\partial_t-\Delta ')a_i$ for all~$i\geq0$. Moreover $a_0(t,x')=\theta(t,x',0)=\sum_{j\geq0} z_j(t)e_j(x')$ by \eqref{D100}
and $a_1=0$ by \eqref{D3}. Therefore for all $i\geq0$,
\begin{IEEEeqnarray*}{rCl}
    a_{2i+1} &=&0\\
    a_{2i} &=& (\partial_t-\Delta')^ia_0\\
    &=& \sum_{j\geq0}(\partial_t-\Delta')^i[z_j(t)e_j(x')]\\
    &=& \sum_{j\geq0}e_j(x')(\partial_t+\lambda_j)^i z_j(t)\\
    &=& \sum_{j\geq0}e_j(x') e^{-\lambda_j t}y_j^{(i)}(t),
\end{IEEEeqnarray*}
where we have set $y_j(t):=e^{\lambda_j t}z_j(t)$. Clearly
\begin{IEEEeqnarray}{rCl}
   \theta(t,x',x_N) &=&\sum_{j\geq0}e^{-\lambda_j t}e_j(x')
   \sum_{i\ge0}y_j^{(i)}(t)\frac{x_N^{2i}}{(2i)!}\label{DD11}
   \\
    u(t,x') &=& \sum_{j\geq0}e^{-\lambda_jt}e_j (x')
    \sum_{i\ge1}\frac{y_j^{(i)}(t)}{(2i-1)!}. \label{DD12}
\end{IEEEeqnarray}
is a formal solution of~\eqref{D1}--\eqref{D3} uniquely defined by the $y_j$'s hence by the~$z_j$'s.

We now give a precise meaning to this formal solution by restricting the~$y_j$'s to be Gevrey of order~$s\in [1,2)$.
\begin{prop}\label{prop2}
Let $s\in [1,2)$, $0< t_1<t_2<\infty$, and consider the sequence $y=(y_j)_{j\ge 0}$ in $C^\infty([0,T])$
satisfying for some constants $M,R >0$
\be
\abs{y_j^{(i)}(t)} \le M \frac{i!^s}{R^i}\qquad\forall i,j\ge 0,\ \forall t\in[t_1,t_2].
\label{C1}
\ee
Then the series in~\eqref{DD11} is uniformly convergent on $[t_1,t_2]\times\overline{\Omega}$ and its sum~$\theta$ is Gevrey of order $s$ in $t$, $1/2$ in $x_1$,\dots,$x_{N-1}$, and $s/2$ in $x_N$; in particular \eqref{D1}--\eqref{D3} and~\eqref{D100} hold for all $(t,x)\in [t_1,t_2]\times\overline{\Omega}$. Moreover the formal control $u$ defined by~\eqref{DD12} is Gevrey of order $s$ in $t$ and $1/2$ in $x_1$,\dots,$x_{N-1}$ on $[t_1,t_2]\times\overline{\omega}$.
\end{prop}
\begin{pf}

We first prove the formal series
\begin{IEEEeqnarray*}{rCl}
    \partial^m_t(\Delta')^{l}\partial^n_{x_N}\theta(t,x',x_N)
    &=&
    \sum_{2i\geq n}\sum_{j\geq0}E_{i,j},
\end{IEEEeqnarray*}
where
\begin{IEEEeqnarray*}{rCl}
    E_{i,j} &:=& \frac{x_N ^{2i-n}}{(2i-n)!} (-\lambda_j)^le_j(x')\partial_t^m \big(e^{-\lambda_j t}y_j^{(i)}(t)\big)
\end{IEEEeqnarray*}
is uniformly convergent on $[t_1,t_2]\times \overline{\Omega}$ with an estimate of the form
\begin{IEEEeqnarray}{rCl}\label{E5}
    \abs{\partial_t^m (\Delta')^l\partial_{x_N}^n\theta(t,x)}
    &\leq&
    C\frac{m!^s}{R_1^m}\,\frac{(2l)!^\frac{1}{2}}{R_1^{2l}}\,\frac{n!^\frac{s}{2}}{R_2^n}
\end{IEEEeqnarray}
for some positive constants $R_1,R_2,R_3$. In what follows $C>0$ denotes a generic constant independent of $m,l,n,t,x$ that may vary from line to line.

Using~\eqref{C1},
\begin{IEEEeqnarray*}{rCl}
    \abs{E_{i,j}} &\leq& \frac{\norm{e_j}_{L^\infty(\omega)}}{(2i-n)!} \sum_{k=0}^m
    \binom{m}{k}\abs{y_j^{(i+k)}(t)}\lambda_j^{l+m-k} e^{-\lambda_jt_1}\\
    &\le& M\frac{m!\norm{e_j}_{L^\infty(\omega)}}{(2i-n)!} \sum_{k=0}^m
  \frac{(i+k)!^s}{k!R^{i+k}} \frac{\lambda_j^{l+m-k} e^{-\lambda_jt_1}}{(m-k)!}.
\end{IEEEeqnarray*}
Pick any $t_0\in(0,t_1)$ such that $\delta:=\frac{t_1-t_0}{2}<2$, and any integer $\kappa>\frac{N-1}{2}$; by the Sobolev inequalities, see e.g.~\cite[Section~7.7]{GilbaT2001book},
\be\label{eq:sobolev}
 \norm{e_j}_{L^\infty(\omega)}
 \le C \norm{e_j}_{H^\kappa(\omega )}
 \le C \lambda_j^\frac{\kappa}{2},
 %\le C j^{\frac{2\kappa}{N-1}}.
\ee
hence
$\norm{e_j}_{L^\infty(\omega)}e^{-\delta\lambda_j} \le C$. It follows that
\begin{IEEEeqnarray*}{rCl}
\frac{\lambda_j^{l+m-k}\norm{e_j}_{L^\infty(\omega )} e^{-2\delta\lambda_j}}{(m - k)!}
&\leq& C\delta^{-(l+m-k)} \frac{(l+m-k)!}{(m-k)!}
\\
&\leq& C\Bigl(\frac{2}{\delta}\Bigr)^{l+m-k}l!,
\end{IEEEeqnarray*}
where we have used $x^pe^{-x}\le p!$  and $(p+q)!\le2^{p+q}p!q!$.
Thus since $s\ge1$
\begin{IEEEeqnarray*}{rCl}
\abs{E_{i,j}} &\le&
Ce^{-\lambda_jt_0}\Bigl(\frac{2}{\delta}\Bigr)^ll!\,  \frac{m!}{(2i-n)!}
\sum_{k=0}^m \frac{2^{m-k+s(i+k)}}{\delta^{m-k}k!}\frac{i!^s k!^s}{R^{i+k}}
\\
&\le& Ce^{-\lambda_jt_0} \Bigl(\frac{2}{\delta}\Bigr)^ll!\,  \frac{i!^s\bigl(\frac{2^s}{R}\bigr)^i}{(2i-n)!} \Bigl(\frac{2}{\delta}\Bigr)^m m!^s \sum_{k=0}^m \Bigl(\frac{2^s}{R}\Bigr)^k .
\end{IEEEeqnarray*}

Now by Stirling's formula
\begin{IEEEeqnarray*}{rClCl}
    l!\Bigl(\frac{2}{\delta}\Bigr)^l &\sim& (2l)!^\frac{1}{2}\frac{(\pi l)^\frac{1}{4}}{\delta^l}
    &\leq& C \frac{(2l)!^\frac{1}{2}}{R_2^{2l}},
\end{IEEEeqnarray*}
where $R_2<\sqrt\delta$. Likewise using also $(p+q)! \le 2^{p+q}p!q!$
\begin{IEEEeqnarray*}{rCl}
    \frac{i!^s\Bigl(\frac{2^s}{R}\Bigr)^i}{(2i-n)!}
    &\sim& \Bigl(\frac{2^s}{R}\Bigr)^i \frac{(2i)!^\frac{s}{2}}{(2i-n)!^\frac{s}{2}}
    \frac{(\pi i)^\frac{s}{4}}{(2i-n)!^{1-\frac{s}{2}}2^{is}}\\
    &\leq& \Bigl(\frac{2^s}{R}\Bigr)^i (\pi i)^\frac{s}{4} \frac{n!^\frac{s}{2}}{(2i-n)!^{1-\frac{s}{2}}}\\
    &\leq& \Bigl(\frac{2^s}{R}\Bigr)^i \Bigl(\frac{n\pi}{2}\Bigr)^\frac{s}{4}
    \biggl(1+\Bigl(\frac{2i-n}{n}\Bigr)^\frac{s}{4}\biggr) \frac{n!^\frac{s}{2}}{(2i-n)!^{1-\frac{s}{2}}}\\
    &\leq& C \frac{n!^\frac{s}{2}}{R_2^n}
    \underbrace{ \frac{1+\Bigl(\frac{2i-n}{n}\Bigr)^\frac{s}{4}}{(2i-n)!^{1-\frac{s}{2}}}
    \sqrt{\frac{2^s}{R}}^{2i-n} }_{A_{2i-n}},
\end{IEEEeqnarray*}
 where $R_3<\sqrt{\frac{R}{2^s}}$ and $\sum_{2i\geq n}{A_{2i-n}}<\infty$. Finally notice
\begin{IEEEeqnarray*}{rCl}
    \Bigl(\frac{2}{\delta}\Bigr)^m \sum_{k=0}^m \Bigl(\frac{2^s}{R}\Bigr)^k
    &\leq& C\frac{1}{R_1^m}
\end{IEEEeqnarray*}
if we set
\begin{IEEEeqnarray*}{rCl}
    R_1 &:=& \begin{cases}
    \frac{\delta R}{2^{s+1}}& \text{if $R<2^s$},\\
    \frac{\delta}{4}& \text{if $R=2^s$},\\
    \frac{\delta}{2}& \text{if $R>2^s$}.
\end{cases}
\end{IEEEeqnarray*}
Collecting the three previous estimates yields
\begin{IEEEeqnarray*}{rCl}
    \sum_{2i\geq n}\sum_{j\geq0}E_{i,j} &\leq&
    C \frac{m!^s}{R_1^m}\,\frac{(2l)!^\frac{1}{2}}{R_2^{2l}}\,\frac{n!^\frac{s}{2}}{R_3^n}
    \sum_{j\geq0}e^{-\lambda_jt_0} \sum_{2i\geq n}{A_{2i-n}}.
\end{IEEEeqnarray*}
The second series in the r.h.s. is convergent and so is the first one by~\eqref{weyl}.
Thus~\eqref{E5} is proved, which shows $\partial_t^m(\Delta')^l\partial_{x_N}^n\theta\in C([t_1,t_2]\times \overline{\Omega})$ for all $m,l,n\ge 0$, hence
$\theta \in C^\infty ([t_1,t_2]\times\overline{\Omega})$ by the Sobolev imbedding theorem.

To complete the proof of the proposition we must extend~\eqref{E5} to every derivative $\partial_t^m\partial_{x'}^{\alpha'}\partial_{x_N}^n\theta$
where $\alpha'=(\alpha_1,\ldots,\alpha _{N-1})\in \N ^{N-1}$ and $\partial_{x'}^{\alpha'}=\partial_{x_1}^{\alpha_1}\cdots\partial_{x_{N-1}}^{\alpha_{N-1}}$.
Pick any integer $p\in(N-1,\infty )$. By the Sobolev inequalities
\be\label{Sob1}
\norm{v}_{L^\infty( \omega )} \le C\norm{v}_{W^{1,p}(\omega)}.
\ee
Repeated  applications of the classical %$L^p$
elliptic estimate
\bes
\norm{v}_{W^{2,p}(\omega)} \leq C\bigl(\norm{\Delta v}_{L^p(\omega)} + \norm{v}_{L^p(\omega)}\bigr)
\ees
(where the constant $C$ depends on $p$) yield for all $l\in\N^\star$
\bes
\norm{v}_{W^{2l,p}(\omega)} \le C \Bigl(\norm{v}_{L^p(\omega)}
+\sum_{k=1}^l\norm{\Delta^kv}_{L^p(\omega)}\Bigr).
\ees
As a consequence and using~\eqref{Sob1},
\be\label{E20}
\norm{v}_{W^{2l,p}(\omega)} \le C \Bigl(\norm{v}_{L^\infty(\omega)}
+\sum_{k=1}^l\norm{\Delta^kv}_{L^\infty(\omega)}\Bigr).
\ee
Pick any $l\in\N^\star$ and any $\alpha'\in\N^{N-1}$ with $\abs{\alpha'}:=\alpha_1+\cdots+\alpha_{N-1}\in\{2l-2,2l-1\}$.
Assume e.g. that $\abs{\alpha'}= 2l-1$, the computations being similar when
$\abs{\alpha'}=2l-2$. Then with \eqref{E20} and~\eqref{E5},
\begin{IEEEeqnarray*}{rCl}
\abs{\partial_t^m\partial_{x'}^{\alpha'}\partial_{x_N}^n\theta(t,x)}
&\leq& \norm{\partial_t^m\partial_{x_N}^n\theta(t,\cdot,x_N)}_{W^{2l-1,\infty}(\omega)}
\\
&\leq&  C\frac{m!^s}{R_1^m}\sum_{k=0}^l\frac{(2k )!^\frac{1}{2}}{R_2^{2k}}\frac{n!^\frac{s}{2}}{R_3^n}
\\
&\leq& C\frac{m!^s}{R_1^m}(l+1)\frac{(2l)!^\frac{1}{2}}{R_2^{2l}}\frac{n!^\frac{s}{2}}{R_3^n}.
\end{IEEEeqnarray*}
Pick $R_2'\in(0,R_2)$. Then $(l+1)\sqrt{2l}\leq C\frac{R_2}{R_2'}^{2l}$ hence
\bes
(l+1)\frac{(2l)!^\frac{1}{2}}{R_2^{2l}}
\le C\frac{(2l-1)!^\frac{1}{2}}{(R_2')^{2l}}
\le C\frac{\alpha_1!^\frac{1}{2}\cdots\alpha_{N-1}!^\frac{1}{2}}
{(\frac{R_2'}{\sqrt{N-1}})^{\alpha_1+\cdots+\alpha_{N-1}}},
\ees
where we have used $(2l-1)!\le(N-1)^{2l-1}\alpha_1!\cdots\alpha_{N-1}!$. Finally, setting $R_2''=\frac{R_2'}{\sqrt{N-1}}$,
\bes
\abs{\partial_t^m\partial_{x'}^{\alpha'}\partial_{x_N}^n\theta(t,x)}
\leq C\frac{m!^s}{R_1^m}
\frac{\alpha_1!^\frac{1}{2}\cdots\alpha_{N-1}!^\frac{1}{2}}{(R_2'')^{\alpha_1+\cdots+\alpha_{N-1}}}
\frac{n!^\frac{s}{2}}{R_3^n},
\ees
which completes the proof of Proposition~\ref{prop2}.\qed \end{pf}

With Proposition~\ref{prop2} at hand, we can derive a constructive null controllability result. The approach is similar to and builds on the proof of Theorem~\ref{thm1}.
\begin{thm}\label{thm3}
Let $\theta_0\in L^2(\Omega)$, $T>0$, $\tau\in(0,T)$ and $s\in(1,2)$. Then there exists a sequence $(y_j)_{j\ge0}$ of functions in
$C^\infty([\tau,T])$ which are Gevrey of order $s$ on~$[\tau, T]$ and such that the control
\bes
u(t,x')=:\begin{cases}
0 & \text{if $0\le t\le\tau$}\\
\sum_{i\ge1,j\ge0}e^{-\lambda_jt}\frac{y_j^{(i)}(t)}{(2i-1)!}e_j(x') &\text{if $\tau <t\le T$}.
\end{cases}
\ees
steers the system~\eqref{D1}--\eqref{D3} from the initial state $\theta(0,.)=\theta_0$ at time~$0$ to the final state~$\theta(T,.)=0$ at time~$T$.

Moreover $u$ is Gevrey of order $s$ in~$t$ and $1/2$ in $x_1,\dots,x_{N-1}$ on $[0,T]\times \overline{\omega}$; %$t\mapsto\theta(t,\cdot)$ is in $C\bigl([0,T],L^2(\Omega)\bigr)$;
the solution $\theta$ of~\eqref{D1}--\eqref{D3} is Gevrey of order~$s$ in~$t$, $1/2$ in $x_1,\dots,x_{N-1}$, and $s/2$ in $x_N$ on~$[\varepsilon,T]\times\overline{\Omega}$ for all $\varepsilon\in(0,T)$.
\end{thm}
\begin{pf}
\paragraph*{\sc Step~1: free evolution}
Let $u(t,x')=0$ for~$(t,x')\in[0,\tau]\times\omega$. Decompose $\theta_0$ as the Fourier series of cosines
\be\label{eq:cjn}
\theta_0(x',x_N)=\sum_{j,n\ge0}c_{j,n}e_j(x')\sqrt{2}\cos(n\pi x_N),
\ee
with convergence in~$L^2(\Omega)$ and $\sum_{j,n\ge0}\abs{c_{j,n}}^2 <\infty$. The solution $\theta$ of~\eqref{D1}--\eqref{D3} starting from~$\theta_0$ then reads
\be\label{WW1}
\theta(t,x',x_N)=\sum_{j,n\ge0}c_{j,n}e^{-(\lambda_j+n^2\pi^2)t}e_j(x')\sqrt{2}\cos(n\pi x_N).
\ee
In particular
\bes%\label{FFF1}
\theta_\tau(t,x):=\theta(\tau,x',x_N)=\sum_{j\ge 0}e^{-\lambda_j\tau}e_j(x')\theta_{\tau,j}(x_N),
\ees
where by Lemma~\ref{lem1}

\begin{IEEEeqnarray*}{rCl}
\theta_{\tau,j}(x_N) &=& \sum_{n\ge0}c_{j,n}e^{-n^2\pi^2\tau}\sqrt{2}\cos(n\pi x_N)
\\
&=& \sum_{i\ge0} \frac{x_N^{2i}}{(2i)!} \underbrace{\left(\sqrt{2}\sum_{n\ge0}c_{j,n} e^{-n^2\pi^2\tau}(-n^2\pi^2)^i\right)}_{=:y_{j,i}}%\IEEEeqnarraynumspace\label{eq:yji}
\end{IEEEeqnarray*}
with $\abs{y_{j,i}} \le C\Bigl(1+\frac{1}{\sqrt\tau}\Bigr)\frac{i!}{\tau^i}$.

\paragraph*{\sc Step~2: construction of the control on~$[\tau,T]$}
By Lemma~\ref{lem2} the functions~$y_j$, $j\ge0$, defined on~$[\tau ,T]$ by
\be\label{eq:yj}
y_j(t):=\phi_s\Bigl(\frac{t-\tau}{T-\tau}\Bigr)\sqrt{2}\sum_{n\ge0}c_{j,n}e^{-n^2\pi^2t}
\ee
are Gevrey of order~$s$ on~$[\tau ,T]$, and satisfy for all $i\ge 0$
\begin{IEEEeqnarray}{rCl}
y_j^{(i)}(\tau) &=& y_{j,i}\label{F2}\\
y_j^{(i)}(T) &=& 0\label{F3}\\
\abs{y_j^{(i)}(t)} &\le&  C\sup_{j,n\ge0}\abs{c_{j,n}}\frac{i!^s}{R^i}, \quad\forall t\in[\tau,T] \label{F1}
\end{IEEEeqnarray}
for some constants $C,R>0$ depending only on~$\tau,T,s$.

The control defined on~$[\tau,T]$ by~\eqref{DD12} then steers the system~\eqref{D1}--\eqref{D3} from $\theta_\tau$ at time~$\tau$ to 0 at time~$T$. Indeed the solution~$\theta$ on~$[\tau,T]$ is given by~\eqref{DD11}, and by~\eqref{F2}-\eqref{F3} satisfies $\theta(\tau,\cdot)=\theta_\tau$ and $\theta(T,\cdot)=0$. By Proposition~\ref{prop2} this solution~$\theta$ is Gevrey of order~$s$ in~$t$, $1/2$ in $x_1,\dots,x_{N-1}$, and $s/2$ in~$x_N$; on the other hand by an easy extension of Lemma~\ref{lem:regstep1} and arguments similar to those in the proof of Proposition~\ref{prop2} the solution $\theta$ defined by~\eqref{WW1} is Gevrey of order~$1$ in~$t$ and $1/2$ in~$x_1,\dots,x_N$ on $[\varepsilon,\tau]\times\overline\Omega$ for all $\varepsilon\in(0,\tau)$.
Since $s>1$, to prove $\theta$ is Gevrey of order $s$ in~$t$, $1/2$ in $x_1,\dots,x_{N-1}$, and $s/2$ in~$x_N$ on $[\varepsilon,T]\times\overline\Omega$ for all $\varepsilon\in(0,T)$ it is then sufficient to check
$\partial_t^k\theta(\tau^+,x)=\partial_t^k\theta(\tau^-,x)$ for $k\geq0$ and~$x\in\overline\Omega$.
But by~\eqref{F2} the series \eqref{DD11} and~\eqref{WW1} coincide at $t=\tau$, hence so do their (space) derivatives. Therefore
\bes
\partial_t^k\theta(\tau^+,x)
=\Delta^{2k}_x\theta(\tau^+,x)
=\Delta^{2k}_x\theta(\tau^-,x)
=\partial_t^k\theta(\tau^-,x).
\ees
%%%%%%%%%%%%%%% OLD EXPLICIT COMPUTATION
%\begin{IEEEeqnarray*}{rCl}
%\IEEEeqnarraymulticol{3}{l}{
%\partial_t^k\theta(\tau^+,x',x_N)
%}\\~~
%& = & \sum_{i,j\ge0} \frac{x_N^{2i}}{(2i)!}e_j(x')
%  \partial_t^k\bigl[e^{-\lambda_jt}y_j^{(i)}(t)\bigr]\Bigl\vert_{t=\tau^+}
%\\
%&=& \sum_{i,j\ge0}\frac{x_N^{2i}}{(2i)!}e_j(x')
%  \sum_{m=0}^k\binom{k}{m}(-\lambda_j)^me^{-\lambda_j\tau}y_{j,i+k-m}
%\\
%&=& \sqrt{2}\sum_{i,j\ge0}\frac{x_N^{2i}}{(2i)!}e_j(x')
%  \sum_{m=0}^k\binom{k}{m}(-\lambda_j)^me^{-\lambda_j\tau}\\&&\quad
%  \sum_{n\ge0}c_{j,n} e^{-n^2\pi^2\tau}(-n^2\pi^2)^i
%\\
%&=& \sqrt{2}\sum_{j,n\ge0}c_{j,n}e^{-(\lambda_j+n^2\pi^2)\tau}e_j(x')\\&&\quad
% \sum_{m=0}^k \binom{k}{m}(-\lambda_j)^m(-n^2\pi^2)^{k-m}
% \sum_{i\ge 0}(-1)^i\frac{(n\pi x_N)^{2i}}{(2i)!}
%\\
%&=& \sum_{j,n\ge0}(-\lambda_j-n^2\pi^2)^kc_{j,n}e^{-(\lambda_j+n^2\pi^2)\tau}e_j(x')\\&&\quad
% \sqrt{2}\cos(n\pi x_N)
%\\
%&=& \partial_t^k\theta(\tau^-,x',x_N).
%\end{IEEEeqnarray*}
%
As a simple consequence $u(t,x')=\partial_{x_N}\theta(t,x',1)$ is Gevrey of order $s$ in~$t$ and $1/2$ in $x_1,\dots,x_{N-1}$ on $[0,T]\times \overline{\omega}$, which completes the proof of Theorem~\ref{thm3}.\qed
\end{pf}

\section{Numerical approximations}\label{sec:estimates}
We now investigate numerical approximations of the control problem~\eqref{D1}--\eqref{D3}. We assume $N\ge2$, the case $N=1$ being similar and simpler. From section~\ref{sec:dimN} we know the control~$u$ achieving the exact control to zero is explicitly given by Theorem~\ref{thm3}, with the functions~$y_j$ and the coefficients~$c_{j,n}$ defined by~\eqref{eq:yj}-\eqref{eq:cjn}; the corresponding solution~$\theta$ is given by~\eqref{WW1} for $t\leq\tau$ and by~\eqref{DD11} for $t\geq\tau$. The aim of this section is to show that approximating the various series by their truncated sums at given $\oi,\oj,\on\in\N$ provides very good approximations, and to give explicit error estimates.

We begin with the free evolution and consider
\bes
\overline\theta(t,x',x_N):=\sqrt{2}\mspace{-10.0mu}\sum_{\substack{0\le j\le\oj\\ 0\le n\le\on}}
c_{j,n}e^{-(\lambda_j+ n^2\pi^2)t}e_j(x')\cos(n\pi x_N)
\ees
instead of the true solution~\eqref{WW1}. Clearly for $0\le t\le \tau$
\bes
\norm{\theta(t)-\overline\theta(t)}^2_{L^2(\Omega)}
=\sum_{(j,n)\not\in[0,\oj]\times[0,\on]}e^{-2(\lambda_j+n^2\pi ^2)t}\abs{c_{j,n}}^2,
\ees
hence $\overline\theta\to\theta$ in~$C([0,\tau ];L^2(\Omega))$ as $\oj,\on\to\infty$. Moreover for $\tau'\in(0,\tau)$ and $r>0$, $\overline\theta\to\theta$ in~$C([\tau',\tau];H^r(\Omega ))$ as $\oj,\on\to\infty$.

The main result of this section concerns the convergence for $t\in[\tau,T]$ of the truncated sums
\begin{IEEEeqnarray*}{rCl}
\overline u(t,x') &:=& \sum_{0\le j\le\oj}e^{-\lambda_jt}e_j(x')\sum_{1\le i\le\oi}
\frac{\overline{y_j}^{(i)}(t)}{(2i-1)!}%\label{eq:ubar}\IEEEeqnarraynumspace
\\
\overline{\theta}(t,x',x_N) &:=& \sum_{0\le j\le\oj}e^{-\lambda_jt}e_j(x')\sum_{0\le i\le\oi} \overline{y_j}^{(i)}(t)\frac{x_N^{2i}}{(2i)!}%\label{eq:thetabar}
\\
\overline{y_j}(t) &:=& \phi_s\Bigl(\frac{t-\tau}{T-\tau}\Bigr)\sqrt2\sum_{0\le n\le\on}c_{j,n}e^{-n^2\pi^2t}. %\label{Z8}
\end{IEEEeqnarray*}
Let also $\hat u$ be the control defined by~$0$ for $t\in[0,\tau]$ and by~$\overline u$ for $t\in(\tau,T]$, and $\hat\theta$ be the solution of~\eqref{D1}--\eqref{D3} when $u$ is replaced by $\hat u$ (still starting from $\theta_0$ at time 0); $\hat\theta$ is thus the solution of the ``true'' control problem when using the approximated control~$\hat u$ instead of the ``ideal'' control $u$. Notice $\hat\theta$ and~$\theta$ agree on $[0,\tau]$, but $\hat\theta$ and $\overline\theta$ differ on $(\tau,T]$ (i.e. the solution produced by the truncated control is not the truncated solution).
\begin{thm}\label{thm4}
Let $\theta_0\in L^2(\Omega)$, $T>0$, $\tau\in(0,T)$, $s\in(1,2)$ and $N\ge2$. Then there are constants $C_1,C_2,C_3,C_4$ such that for all $t\in[\tau,T]$
\be\label{ZZ1}
\norm{\theta(t)-\overline\theta(t)}_{L^\infty(\Omega)} \le C_1f(\oi,\oj,\on)\norm{\theta_0}_{L^2(\Omega)},
\ee
where
\be\label{eq:Deltas}
f(\oi,\oj,\on):=e^{-C_2\,\oj^\frac{2}{N-1}} + e^{-C_3\,\oi\ln\oi} + e^{-C_4\,\on^2}.
\ee
In~\eqref{ZZ1}-\eqref{eq:Deltas} one may choose $C_2<A_1\tau$ (see~\eqref{weyl}), $C_3<2-s$, and $C_4<\pi^2\tau$, while
$C_1=C_1(N,\omega,\tau,s,C_2,C_3,C_4)$.
\end{thm}
\begin{pf} For $(t,x',x_N)\in[\tau ,T]\times\omega\times[0,1]$
\be\label{delta}
\abs{\theta(t,x',x_N)-\overline\theta(t,x',x_N)} \le \Delta_1+\Delta_2+\Delta_3,
\ee
where
\begin{IEEEeqnarray*}{rCl}
\Delta_1 &:=&
\sum_{j>\oj}e^{-\lambda_jt}\abs{e_j(x')}\sum_{i\ge0}\frac{\bigl\vert y_j^{(i)}(t)\bigr\vert}{(2i)!}
\\
\Delta_2 &:=&
\sum_{0\le j\le\oj}e^{-\lambda_jt}\abs{e_j(x')}\sum_{i>\oi}\frac{\bigl\vert y_j^{(i)}(t)\bigr\vert}{(2i)!}
\\
\Delta _3 &:=&
\sum_{0\le j\le\oj}e^{-\lambda_jt}\abs{e_j(x')}
\sum_{0\le i\le\oi}
\frac{\Bigl\vert\partial_t^i\bigl[\phi(t)\sum\limits_{n>\on}c_{j,n}e^{-n^2\pi^2t}\bigr]\Bigr\vert}{(2i)!}
\end{IEEEeqnarray*}
and $\phi(t):=\sqrt2~\phi_s\bigl(\frac{t-\tau}{T-\tau}\bigr)$.

We first bound~$\Delta_1$. On the one hand by~\eqref{F1} and $s<2$
\begin{IEEEeqnarray}{rCl}
\sum_{i\ge0}\frac{\abs{y_j^{(i)}(t)}}{(2i)!}
&\le& C\sup_{j,n\ge0}\abs{c_{j,n}}\sum_{i\ge0}\frac{(i!)^s}{(2i)!R^i}\nonumber\\
&\le& C\sup_{j,n\ge0}\abs{c_{j,n}}\label{Z23}
\end{IEEEeqnarray}
for some constant $C=C(\tau,T,s)$. On the other hand using \eqref{eq:sobolev}, \eqref{weyl} and $x^pe^{-x}\le p!$ and taking $\rho\in(0,1)$
\begin{IEEEeqnarray}{rCl}
e^{-\lambda_j\tau}\norm{e_j}_{L^\infty(\omega)}
&\le& Ce^{-\lambda_j\tau}\lambda_j^\frac{\kappa}{2}\nonumber\\
&\le& Ce^{-\lambda_j\rho\tau}\nonumber\\
&\le& Ce^{-A_1\rho\tau j^\frac{2}{N-1}},\label{Z24}
\end{IEEEeqnarray}
where the constant $C$ eventually depends on~$\tau,N$. But
\begin{IEEEeqnarray*}{rCl}
\sum_{j>\oj}e^{-Kj^r}
&\le& \int_{\oj}^{\infty}e^{-Kx^r}dx
\\
&=& \frac{1}{r}\int_{\oj^r}^\infty e^{-Ky}y^{\frac{1}{r}-1}dy
\\
&\le& \begin{cases}
    \frac{\oj^{1-r}}{rK}e^{-K\oj^r}& \text{if $r\ge1$}\\
    K'e^{-K''\oj^r}& \text{if $r<1$},
\end{cases}%\label{Z25}
\end{IEEEeqnarray*}
where $r,K,K'>0$ and $K''\in (0,K)$. Picking $C_2<A_1\tau$ then yields
\be
\label{Z28}
\sum_{j>\oj } e^{-\lambda _j \tau } ||e_j||_{L^\infty ( \omega ) }
\le C e^{-C_2 \oj ^\frac{2}{N-1}},
\ee
where $C= C(\tau,N,C_2)$. We conclude by \eqref{Z23} and~\eqref{Z28}
\begin{IEEEeqnarray}{rCl}
\Delta_1
&\le& \sum_{j>\oj}e^{-\lambda_j \tau}\norm{e_j}_{L^\infty(\omega)}
\sum_{i\ge0}\frac{\abs{y_j^{(i)}(t)}}{(2i)!}\nonumber
\\
&\le& Ce^{-C_2\oj^\frac{2}{N-1}}\sup_{j,n\ge0}\abs{c_{j,n}},\label{delta1}
\end{IEEEeqnarray}
where $C=C(\tau,T,N,C_2,s)$.

We next bound~$\Delta_2$. By~\eqref{F1} and using Stirling's formula for $i!$ and~$(2i)!$
\begin{IEEEeqnarray*}{rCl}
\sum_{i\ge\oi}\frac{\abs{y_j^{(i)}(t)}}{(2i)!}
&\le& C\sup_{j,n\ge0}\abs{c_{j,n}}\sum_{i\ge\oi}\frac{(i!)^s}{(2i)!R^i}
\\
&\le& C\sup_{j,n\ge0}\abs{c_{j,n}}
\sum_{i\ge\oi}\frac{i^\frac{s-1}{2}}{\left(\frac{i}{e}\right)^{(2-s)i}(4R)^i}
\\
&\le& C\sup_{j,n\ge0}\abs{c_{j,n}}
\sum_{i\ge\oi}\frac{i^\frac{s-1}{2}e^{-(2-s)i(\ln i-1)}}{(4R)^i}
\end{IEEEeqnarray*}
for some constant $C=C(\tau,T,s)$. Pick $C_3<2-s$ and $\sigma\in (C_3, 2-s)$; then
\bes
\frac{i^\frac{s-1}{2}e^{(2-s)(1-\ln i)i}}{(4R)^i}
%=\biggl(\frac{i^\frac{s-1}{2}}{(4R)^ie^{(2-s-\sigma)i(\ln i-1)}}\biggr)e^{-\sigma i(\ln i-1)}
\le C e^{-\sigma i(\ln i-1)}
\ees
where $C=C(s,\sigma,R)$, and
\begin{IEEEeqnarray*}{rCl}
\sum_{i>\oi}e^{-\sigma i(\ln i-1)}
&\le& \int_{\oi}^\infty e^{-\sigma x(\ln x-1)}dx
\\
&\le& C\int_{\oi(\ln\oi-1)}^\infty e^{-\sigma y}dy
\\
&\le& Ce^{-C_3\oi\ln\oi}.
\end{IEEEeqnarray*}
We conclude by~\eqref{Z24}
\begin{IEEEeqnarray}{rCl}
\Delta_2
&\le& \sum_{0\le j\le\oj}e^{-\lambda_j\tau}\norm{e_j}_{L^\infty(\omega)}
\sum_{i\ge\oi}\frac{\bigl\vert y_j^{(i)}(t)\bigr\vert}{(2i)!}\nonumber
\\
&\le& C\sup_{j,n\ge0}\abs{c_{j,n}}e^{-C_3\oi\ln\oi}
\sum_{j\ge0}e^{-A_1\rho\tau j^\frac{2}{N-1}}\nonumber
\\
&\leq& C\sup_{j,n\ge0}\abs{c_{j,n}}e^{-C_3\oi\ln\oi}\label{delta2},
\end{IEEEeqnarray}
where $C=C(\tau,T,N,s,C_3)$, since the last series clearly converges.

We finally bound~$\Delta_3$. Pick $\rho\in(0,1)$ and notice first
\bes
\sum_{n>\on}e^{-n^2\pi^2\rho\tau}
\le e^{-\on^2\pi^2\rho\tau}\sum_{p\ge1}e^{-p^2\pi^2\rho\tau}
\le Ce^{-\on ^2\pi^2\rho\tau},
\ees
the series being clearly convergent. Then for $t\ge\tau$
\begin{IEEEeqnarray*}{rCl}
\Bigl\vert\partial_t^i\sum_{n>\on}c_{j,n}e^{-n^2\pi^2t}\Bigr\vert
&=& \Bigl\vert\sum_{n>\on}c_{j,n}(-n^2\pi^2)^ie^{-n^2\pi^2t}\Bigr\vert
\\
&\le& \sup_{j,n\ge0}\abs{c_{j,n}}\sum_{n>\on}(n^2\pi^2)^ie^{-n^2\pi^2\tau}
\\
&\le& i!\frac{\sup_{j,n\ge0}\abs{c_{j,n}}}{\bigl((1-\rho)\tau\bigr)^i}\sum_{n>\on}e^{-n^2\pi^2\rho\tau}
\\
&\le& C\sup_{j,n\ge0}\abs{c_{j,n}}\frac{i!e^{-\on ^2\pi^2\rho\tau}}{\bigl((1-\rho)\tau\bigr)^i},
\end{IEEEeqnarray*}
where we have used the estimate $x^i/i!\le e^x$ for $x=n^2\pi^2(1-\rho)\tau$. In other words $t\mapsto\sum_{n>\on}c_{j,n}e^{-n^2\pi^2t}$ is Gevrey of order~$1$ hence~$s$ on~$[\tau,T]$. Since the product of functions Gevrey of order~$s$ is also Gevrey of order~$s$
\bes
\Bigl\vert\partial_t^i\Bigl[\phi(t)\sum_{n>\on}c_{j,n}e^{-n^2\pi^2t}\Bigr]\Bigr\vert
\le Ce^{-\on^2\pi^2\rho\tau}\sup_{j,n\ge0}\abs{c_{j,n}}\frac{i!^s}{R^i}
\ees
for some $C,R$ depending on $\tau,T,\rho,s$. We conclude by~\eqref{Z24}
\begin{IEEEeqnarray}{rCl}
\Delta_3
&\le& \sum_{0\le j\le\oj}e^{-\lambda_j\tau }\norm{e_j}_{L^\infty(\omega)}
\sum_{0\le i\le\oi}C\sup_{j,n\ge0}\abs{c_{j,n}}\frac{i!^se^{-C_4\on^2}}{(2i)!R^i}\nonumber
\\
&\le& Ce^{-C_4\on^2}\sup_{j,n\ge0}\abs{c_{j,n}}
\sum_{j\ge0}e^{-A_1\rho\tau j^\frac{2}{N-1}}\sum_{i\ge0}\frac{i!^s}{(2i)!R^i}\nonumber
\\
&\leq& Ce^{-C_4\on^2}\sup_{j,n\ge0}\abs{c_{j,n}},\label{delta3}
\end{IEEEeqnarray}
where $C=C(\tau,T,\rho,s)$, since the last two series clearly converge; we have set $C_4:=\pi^2\rho\tau$.

Collecting \eqref{delta}, \eqref{delta1}, \eqref{delta2} and~\eqref{delta3} then yields~\eqref{ZZ1}.\qed
\end{pf}
\begin{cor}\label{cor4}
Using the notations of Theorem~\ref{thm4} and $q\in\N^\star$, there are constants $C_1',C_1''$ depending on $N,\omega,\tau,s,C_2,C_3,C_4,q$ such that
\begin{IEEEeqnarray}{rCl}
\norm{\theta-\overline\theta}_{W^{q,\infty}((\tau,T)\times\Omega)}
&\le& C_1'f(\oi,\oj,\on)\norm{\theta_0}_{L^2(\Omega)}\label{FF1}
\\
\norm{u-\overline u}_{W^{q-1,\infty}((\tau,T)\times\Omega)}
&\le& C_1''f(\oi,\oj,\on)\norm{\theta_0}_{L^2(\Omega)}\IEEEeqnarraynumspace\label{FF2}.
\end{IEEEeqnarray}
\end{cor}
\begin{pf}
The proof of~\eqref{FF1} runs along the same lines as the proof of Theorem~\ref{thm4}. Indeed any additional factor in the $\Delta_i$'s produced by a derivation in $t$, $x'$ or $x_N$ is absorbed by $e^{-\lambda_j\tau}$ or $i!^{s-2}$. On the other hand \eqref{FF2} follows at once from \eqref{FF1} and the fact that $\overline u=\partial_{x_N}\overline\theta(t,x',1)$.
\end{pf}
\begin{cor}\label{cor5}
Using the notations of Theorem~\ref{thm4}, there is a constant $C_1'''=C_1'''(N,\omega,\tau,s,C_2,C_3,C_4)$ such that
\begin{IEEEeqnarray}{rCl}
\lVert\theta-\hat\theta\rVert_{L^{\infty}((0,T)\times\Omega)}
&\le& C_1'''f(\oi,\oj,\on)\norm{\theta_0}_{L^2(\Omega)}\label{FF4}
\end{IEEEeqnarray}
\end{cor}
\begin{pf}
Let
\bes
\tilde\theta:=\theta-\hat\theta-\frac{x_N^2}{2}\bigl(u(t,x')-\overline{u}(t,x')\bigr);
\ees
$\tilde\theta$ is the solution of the initial–boundary-value problem
\begin{IEEEeqnarray*}{rCl'l}
    \IEEEeqnarraymulticol{3}{l}{
    \partial_t\theta(t,x)-\Delta\tilde\theta(t,x)
    }\nonumber\\
    &=&(\Delta-\partial_t)\bigl[\frac{x_N^2}{2}(u-\overline{u})\bigr], &t\in(0,T),x\in\Omega %\label{DD1}
    \\
    \partial_\nu\tilde\theta(t,x) &=& 0, &t\in(0,T),x\in\partial\Omega%\label{DD3}
    \\
    \tilde\theta(0,x) &=&0, &x\in\Omega.%\label{DD4}
\end{IEEEeqnarray*}
Then \eqref{FF4} follows from~\eqref{FF2} with $q=3$ and classical semigroup estimates in~$C(\overline{\Omega})$.
\end{pf}

\section{Numerical experiments}\label{sec:numerics}
We illustrate the approach on two numerical examples, in dimension 1 and~2.

\subsection{Example in dimension~1}
We have taken $T=0.35$, $\tau=0.05$, $s=1.65$. The initial condition $\theta_0$ is a step function with $\theta_0(x)=-0.75$ on~$[0,1/2)$ and $\theta_0(x)=1.25$ on~$(1/2,1)$, with nonzero Fourier coefficients $c_0=\frac{\sqrt2}{8}$ and $c_{2p+1}=\frac{(-1)^{p+1}}{2p+1}\frac{2\sqrt2}{\pi}$ for $p\geq0$; notice $\theta_0$ is not continuous and that its Fourier coefficients decay fairly slowly. The series for the control~$u$ and the function~$y$ of Theorem~\ref{thm1} have been truncated at a ``large enough'' order for a good accuracy, namely $\oi=35$ and~$\on=25$.
For the Gevrey ``step function'' we have used
\bes
\phi_s(t):=\begin{cases}
1 & \text{if $t\leq0$}\\
0 & \text{if $t\geq1$}\\
1-\dfrac{\int_0^t\varphi_s(\rho)d\rho}{\int_0^1\varphi_s(\rho)d\rho}
&\text{if $t\in(0,1)$},
\end{cases}
\ees
where $k = (s-1)^{-1}$ and $\varphi_s$ is the ``bump function''
 \bes
\varphi_s(t):=\begin{cases}
0 & \text{if $t\not\in(0,1)$}\\
\exp\Bigl(\frac{-1}{2t^k(1-t)^k}\Bigr)
&\text{if $t\in(0,1)$}.
\end{cases}
\ees
This ``step function'' was preferred to the one in section~\ref{eq:preliminaries} because its derivatives are easier to compute formally.

\begin{figure}
%\begin{center}
\hspace{-9mm}\includegraphics[width=1.2\columnwidth]{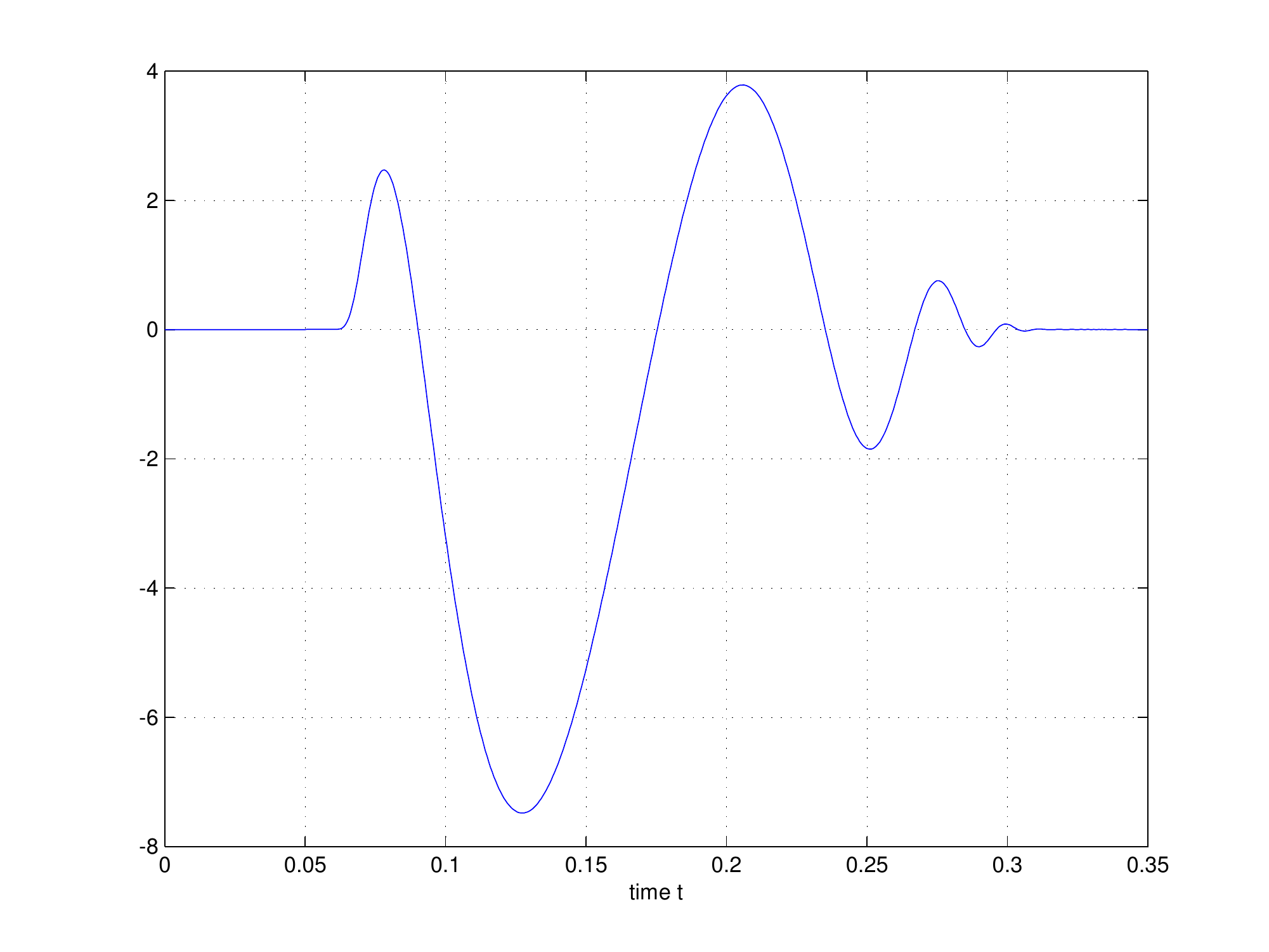}
\caption{1D case: time evolution of the control~$u(t)$.}
\label{fig:u1D}
%\end{center}
\end{figure}

\begin{figure}
%\begin{center}
\hspace{-9mm}\includegraphics[width=1.2\columnwidth]{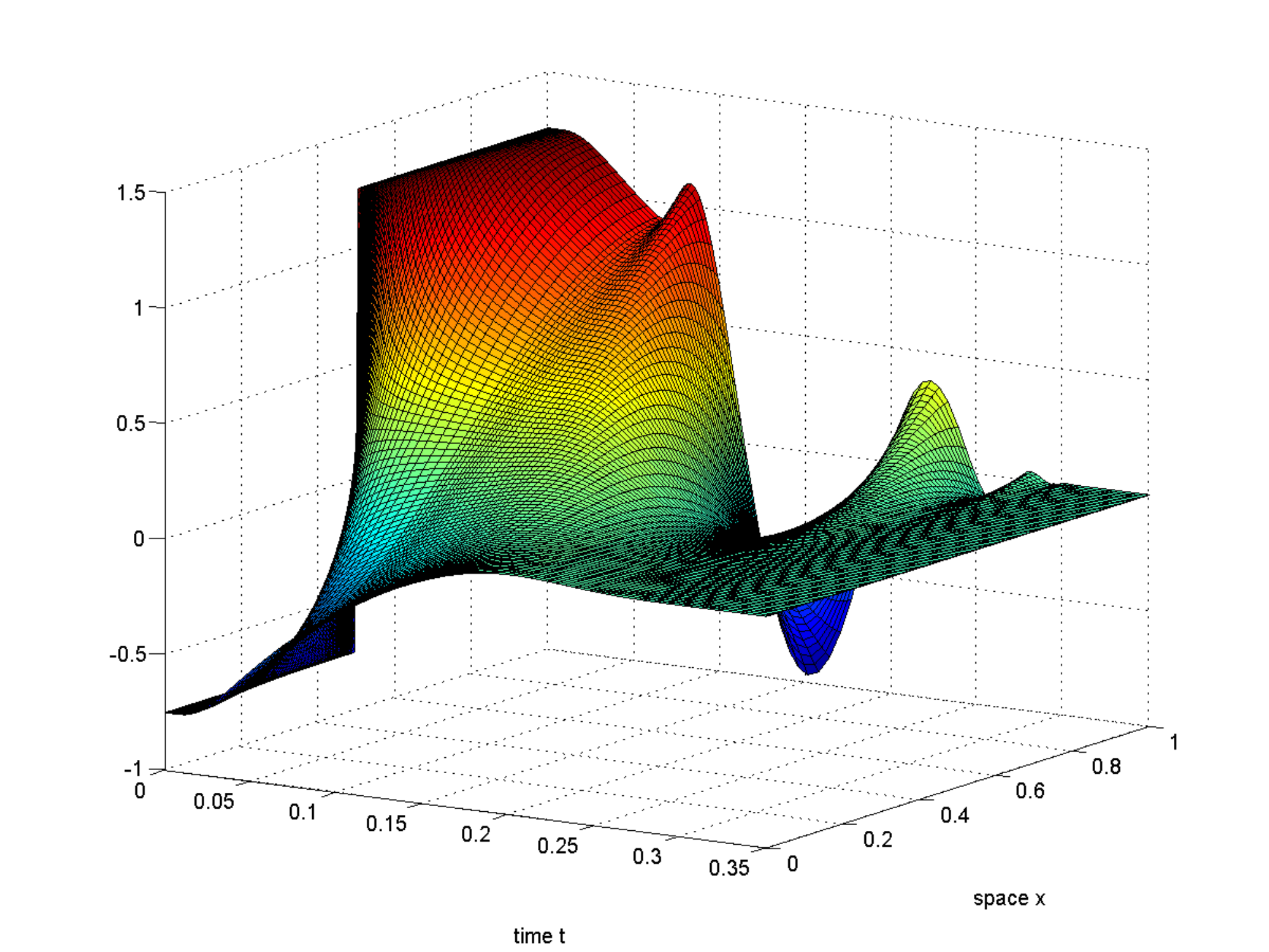}
\caption{1D case: time evolution of the temperature~$\theta(t,x)$.}
\label{fig:theta1D}
%\end{center}
\end{figure}

Fig.~\ref{fig:u1D} shows the (truncated) control $\hat u$ given by Theorem~\ref{thm1}, and Fig.~\ref{fig:theta1D} the resulting simulated temperature~$\hat\theta$ (using the notations of section~\ref{sec:estimates}). The simulation consists of a finite-difference semi-discretization in space with 100 cells, the resulting set of ODEs being integrated with a stiff implicit solver (Simulink ode15s); as expected $\hat\theta$ closely agrees with the truncated series~$\overline\theta$.

\subsection{Example in dimension~2}
\begin{figure}
%\begin{center}
\hspace{-9mm}\includegraphics[width=1.2\columnwidth]{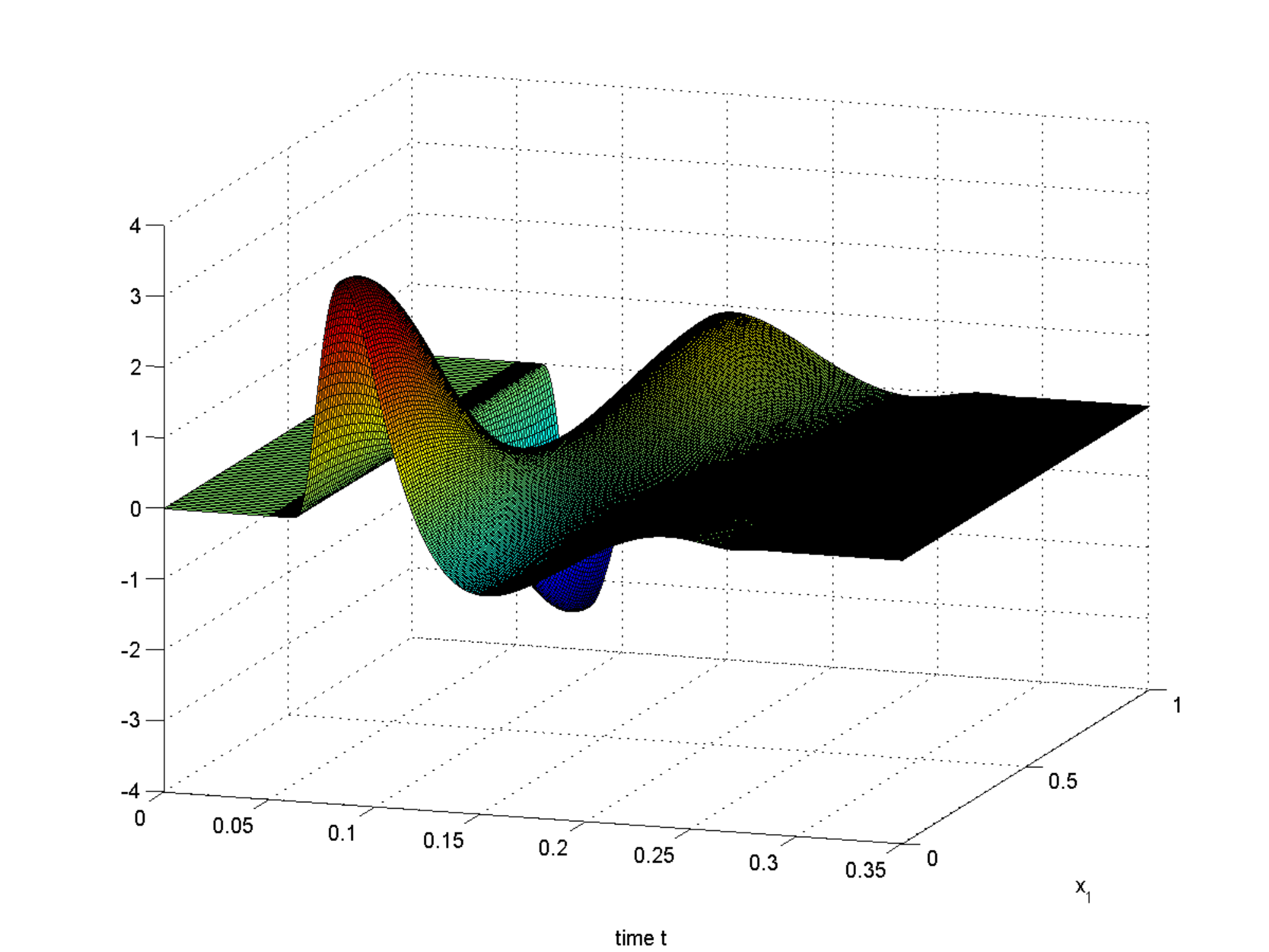}
\caption{2D case: time evolution of the control~$u(t,x_1)$.}
\label{fig:u2D}
%\end{center}
\end{figure}

$T,\tau,s,\oi,\on$ are the same as previously, with moreover $\oj=25$ and $L=1$. The initial condition is the ``double step''
\begin{IEEEeqnarray*}{rCl}
    \theta_0(x_1,x_2) &:=& \begin{cases}
    -1& \text{if $(x_1,x_2)\in(0,\frac{1}{2})\times(0,\frac{1}{2})$},\\
    1& \text{if $(x_1,x_2)\in(0,\frac{1}{2})\times(\frac{1}{2},1)$},\\
    1& \text{if $(x_1,x_2)\in(\frac{1}{2},1)\times(0,\frac{1}{2})$},\\
    -1& \text{if $(x_1,x_2)\in(\frac{1}{2},1)\times(\frac{1}{2},1)$}.
\end{cases}
\end{IEEEeqnarray*}
Its nonzero Fourier coefficients are
\begin{IEEEeqnarray*}{rCl'l}
    c_{2l+1,2p+1} &=& -2\frac{(-1)^{l+p}}{\pi^2(2l+1)(2p+1)}, & l,p\geq0.
\end{IEEEeqnarray*}
Fig.~\ref{fig:u2D} shows the control~$u(t,x_1)$ given by Theorem~\ref{thm3}.

%===============================================================================

\bibliographystyle{abbrv}        % Include this if you use bibtex
\bibliography{automatica}           % and a bib file to produce the

\begin{thebibliography}{10}

\bibitem{BelgaK2011IP}
F.~Belgacem and S.~Kaber.
\newblock On the dirichlet boundary controllability of the one-dimensional heat
  equation: Semi-analytical calculations and ill-posedness degree.
\newblock {\em Inverse Problems}, 27(5):055012, 2011.

\bibitem{BoyerHL2011NM}
F.~Boyer, F.~Hubert, and J.~Le~Rousseau.
\newblock Uniform controllability properties for space/time-discretized
  parabolic equations.
\newblock {\em Numer. Math.}, 118(4):601--661, 2011.

\bibitem{CarthGL1994JOTA}
C.~Carthel, R.~Glowinski, and J.~Lions.
\newblock On exact and approximate boundary controllabilities for the heat
  equation: A numerical approach.
\newblock {\em Journal of Optimization Theory and Applications},
  82(3):429--484, 1994.

\bibitem{FattoR1971ARMA}
H.~Fattorini and D.~Russell.
\newblock Exact controllability theorems for linear parabolic equations in one
  space dimension.
\newblock {\em Arch. Rational Mech. Anal.}, 43(4):272--292, 1971.

\bibitem{FernaM2011HAL}
E.~Fernández-Cara and A.~Münch.
\newblock {Numerical null controllability of the 1D heat equation: primal
  methods}.
\newblock HAL preprint http://hal.archives-ouvertes.fr/hal-00687884, 2011.

\bibitem{FliesLMR1995IJoC}
M.~Fliess, J.~L\'evine, P.~Martin, and P.~Rouchon.
\newblock Flatness and defect of non-linear systems: Introductory theory and
  examples.
\newblock {\em International Journal of Control}, 61(6):1327--1361, 1995.

\bibitem{FursiI1996book}
A.~V. Fursikov and O.~Y. Imanuvilov.
\newblock {\em Controllability of evolution equations}, volume~34 of {\em
  Lecture Notes Series}.
\newblock Seoul National University Research Institute of Mathematics Global
  Analysis Research Center, 1996.

\bibitem{GarciOTToappearJoIaIP}
G.~Garcia, A.~Osses, and M.~Tapia.
\newblock A heat source reconstruction formula from single internal
  measurements using a family of null controls.
\newblock {\em Journal of Inverse and III-posed Problems}, to appear.

\bibitem{GilbaT2001book}
D.~Gilbarg and N.~Trudinger.
\newblock {\em Elliptic partial differential equations of second order}.
\newblock Classics in Mathematics. Springer-Verlag, Berlin, 2001.
\newblock Reprint of the 1998 edition.

\bibitem{Horma1983book}
L.~H{\"o}rmander.
\newblock {\em The analysis of linear partial differential operators. {I}}.
\newblock Springer-Verlag, 1983.

\bibitem{Jones1977JMAA}
B.~Jones~Jr.
\newblock A fundamental solution for the heat equation which is supported in a
  strip.
\newblock {\em J. Math. Anal. Appl.}, 60(2):314--324, 1977.

\bibitem{LarocMR2000IJRNC}
B.~Laroche, P.~Martin, and P.~Rouchon.
\newblock Motion planning for the heat equation.
\newblock {\em Int J Robust Nonlinear Control}, 10(8):629--643, 2000.

\bibitem{LebeaR1995CPDE}
G.~Lebeau and L.~Robbiano.
\newblock Contr\^ole exact de l'\'equation de la chaleur.
\newblock {\em Comm. Partial Differential Equations}, 20(1-2):335--356, 1995.

\bibitem{LinL1995AMO}
Y.-J. Lin~Guo and W.~Littman.
\newblock Null boundary controllability for semilinear heat equations.
\newblock {\em Appl Math Optim}, 32(3):281--316, 1995.

\bibitem{Littm1978ASNSPCS}
W.~Littman.
\newblock Boundary control theory for hyperbolic and parabolic partial
  differential equations with constant coefficients.
\newblock {\em Ann. Scuola Norm. Sup. Pisa Cl. Sci. (4)}, 5(3):567--580, 1978.

\bibitem{LittmT2007}
W.~Littman and S.~Taylor.
\newblock The heat and {S}chr\"odinger equations: boundary control with one
  shot.
\newblock In {\em Control methods in {PDE}-dynamical systems}, volume 426 of
  {\em Contemp. Math.}, pages 293--305. Amer. Math. Soc., Providence, RI, 2007.

\bibitem{LuxemK1971TAMS}
W.~A.~J. Luxemburg and J.~Korevaar.
\newblock Entire functions and {M}\"untz-{S}z\'asz type approximation.
\newblock {\em Trans. Amer. Math. Soc.}, 157:23--37, 1971.

\bibitem{LynchR2002IJC}
A.~Lynch and J.~Rudolph.
\newblock Flatness-based boundary control of a class of quasilinear parabolic
  distributed parameter systems.
\newblock {\em Int J Control}, 75(15):1219--1230, 2002.

\bibitem{MartiRR2013CPDE}
P.~Martin, L.~Rosier, and P.~Rouchon.
\newblock Null controllability of the 1{D} heat equation using flatness.
\newblock In {\em 1st IFAC workshop on Control of Systems Governed by Partial
  Differential Equations (CPDE2013)}, 2013.

\bibitem{MartiRR2013CDC}
P.~Martin, L.~Rosier, and P.~Rouchon.
\newblock Null controllability of the 2{D} heat equation using flatness.
\newblock In {\em 52nd IEEE Conference on Decision and Control}, 2013.

\bibitem{Meure2011A}
T.~Meurer.
\newblock Flatness-based trajectory planning for diffusion-reaction systems in
  a parallelepipedona spectral approach.
\newblock {\em Automatica}, 47(5):935--949, 2011.

\bibitem{MeureZ2008MCMDS}
T.~Meurer and M.~Zeitz.
\newblock Model inversion of boundary controlled parabolic partial differential
  equations using summability methods.
\newblock {\em Math. Comput. Model. Dyn. Syst.}, 14(3):213--230, 2008.

\bibitem{MicuZ2011SCL}
S.~Micu and E.~Zuazua.
\newblock Regularity issues for the null-controllability of the linear 1-d heat
  equation.
\newblock {\em Syst Control Lett}, 60(6):406--413, 2011.

\bibitem{MunchZ2010IP}
A.~M{\"u}nch and E.~Zuazua.
\newblock Numerical approximation of null controls for the heat equation:
  ill-posedness and remedies.
\newblock {\em Inverse Problems}, 26(8):085018, 39, 2010.

\bibitem{Ramis1978}
J.-P. Ramis.
\newblock D\'evissage {G}evrey.
\newblock In {\em Journ\'ees {S}inguli\`eres de {D}ijon ({U}niv. {D}ijon,
  {D}ijon, 1978)}, volume~59 of {\em Ast\'erisque}, pages 4, 173--204. Soc.
  Math. France, Paris, 1978.

\bibitem{Roe1998book}
J.~Roe.
\newblock {\em Elliptic operators, topology and asymptotic methods}, volume 395
  of {\em Pitman Research Notes in Mathematics Series}.
\newblock Longman, Harlow, second edition, 1998.

\bibitem{Rosie2002CAM}
L.~Rosier.
\newblock A fundamental solution supported in a strip for a dispersive
  equation.
\newblock {\em Comput. Appl. Math.}, 21(1):355--367, 2002.
\newblock Special issue in memory of Jacques-Louis Lions.

\bibitem{Rudin1987book}
W.~Rudin.
\newblock {\em Real and complex analysis}.
\newblock McGraw-Hill Book Co., third edition, 1987.

\bibitem{Shubi2001book}
M.~A. Shubin.
\newblock {\em Pseudodifferential operators and spectral theory}.
\newblock Springer-Verlag, Berlin, second edition, 2001.

\bibitem{Yaman1989AGAG}
T.~Yamanaka.
\newblock A new higher order chain rule and gevrey class.
\newblock {\em Annals of Global Analysis and Geometry}, 7(3):179--203, 1989.

\bibitem{Zheng2008AA}
C.~Zheng.
\newblock Controllability of the time discrete heat equation.
\newblock {\em Asymptotic Anal}, 59(3-4):139--177, 2008.

\bibitem{Zuazu2006}
E.~Zuazua.
\newblock Control and numerical approximation of the wave and heat equations.
\newblock In {\em International {C}ongress of {M}athematicians. {V}ol. {III}},
  pages 1389--1417. Eur. Math. Soc., Z\"urich, 2006.

\end{thebibliography}
                                 % bibliography (preferred). The
                                 % correct style is generated by
                                 % Elsevier at the time of printing.

%\appendix
%\section{A summary of Latin grammar}    % Each appendix must have a short title.
%\section{Some Latin vocabulary}         % Sections and subsections are supported
                                        % in the appendices.
\end{document}